\documentclass[journal,twocolumn,dvipdfm]{IEEEtran}
\usepackage{amsmath,amssymb,epsf,epsfig,times}
\usepackage[all]{xy}
\usepackage{graphicx,color}
\usepackage{subfigure}
\usepackage{url}
\usepackage{cite}
\usepackage{keyval,times}
\usepackage[dvipdfm]{hyperref}
\hypersetup{%
   pdftitle={Cooperative Control},
   pdfauthor={Y. Cao, W. Yu, W. Ren and G. Chen},
   pdfkeywords={Cooperative Control, information consensus,
      multi-agent system.},
   bookmarksnumbered,
   pdfstartview={FitH},
   colorlinks=true,
   breaklinks=true,
   citecolor=blue,
}%

\newcommand{\re}{\mathbb{R}}

\newcommand{\norm}[1]{\left\|#1\right\|}

\newcommand{\Acal}{\mathcal{A}}

\newcommand{\Ucal}{\mathcal{U}}

\newcommand{\Lcal}{\mathcal{L}}

\newcommand{\Gcal}{\mathcal{G}}

\newcommand{\Zcal}{\mathcal{Z}}

\newcommand{\onebf}{\mathbf{1}}
\newcommand{\zerobf}{\mathbf{0}}

\newcommand{\ith}{i{\text{th}}}

\newcommand{\OMIT}[1]{}
\begin{document}
\title{An Overview of Recent Progress in the Study of Distributed Multi-agent Coordination}
\author{Yongcan~Cao,~\IEEEmembership{Member,~IEEE,} Wenwu Yu,~\IEEEmembership{Member,~IEEE,}\\
        Wei~Ren,~\IEEEmembership{Member,~IEEE}, and Guanrong Chen,~\IEEEmembership{Fellow,~IEEE}
\thanks{Y. Cao is with the Control Science Center of Excellence, Air Force Research Laboratory, Wright-Patterson AFB, OH 45433, USA. W. Yu is with the Department of Mathematics, Southeast University, Nanjing 210096, China and also with the School of Electrical and Computer Engineering, RMIT
University, Melbourne VIC 3001, Australia. W.~Ren is with the Department of Electrical Engineering, University of California, Riverside, CA 92521, USA. G. Chen is
with the Department of Electronic Engineering, City University of Hong Kong, Hong Kong SAR, China.}
}


\maketitle

\begin{abstract}
This article reviews some main results and progress in distributed multi-agent coordination, focusing on papers published in major control systems and robotics journals since 2006. Distributed coordination of multiple vehicles, including unmanned aerial vehicles, unmanned ground vehicles and unmanned underwater vehicles, has been a very active research subject studied extensively by the systems and control community. The recent results in this area are categorized into several directions, such as consensus, formation control, optimization, task assignment, and estimation. After the review, a short discussion section is included to summarize the existing research and to propose several promising research directions along with some open problems that are deemed important for further investigations.
\end{abstract}


\IEEEpeerreviewmaketitle

\section{Introduction}

\IEEEPARstart{C}{ontrol} theory and practice may date back to the beginning of the last century when Wright Brothers attempted their first test flight in 1903. Since then, control theory has gradually gained popularity, receiving more and wider attention especially during the World War II when it was developed and
applied to fire-control systems, missile navigation and guidance, as well as various electronic automation devices. In the past several decades, modern control theory was further advanced due to the booming of aerospace technology based on large-scale engineering systems.

During the rapid and sustained development of the modern control theory, technology for controlling a single vehicle, albeit higher-dimensional and complex, has become relatively mature and has produced many effective tools such as PID control, adaptive control, nonlinear control, intelligent control, and robust control methodologies. In the past two decades in particular, control of multiple vehicles has received increasing demands spurred by the fact that many benefits can be obtained when a single complicated vehicle is equivalently replaced by multiple yet simpler vehicles. In this endeavor, two approaches are commonly adopted for controlling multiple vehicles: a
centralized approach and a distributed approach. The centralized approach is based on the assumption that a central station is available and powerful enough to control a whole group of vehicles. Essentially, the centralized approach is a direct extension of the traditional single-vehicle-based control philosophy and strategy. On the contrary, the distributed approach does not require a central station for control, at the cost of becoming far more complex in structure and organization. Although both approaches are considered practical depending on the situations and conditions of the real applications, the distributed approach is believed more promising due to many inevitable physical constraints such as limited resources and energy, short wireless communication ranges, narrow bandwidths, and large sizes of vehicles to manage and control. Therefore, the focus of this overview is placed on the distributed approach.

In distributed control of a group of autonomous vehicles, the main objective typically is to have the
whole group of vehicles working in a cooperative fashion throughout a distributed protocol. Here, {\it cooperative\/} refers to a close relationship among all vehicles in the group where {\it information sharing\/} plays a central role. The distributed approach has many advantages in achieving cooperative group performances, especially with low operational costs, less system requirements, high robustness, strong adaptivity, and flexible scalability, therefore has been widely recognized and appreciated.

The study of distributed control of multiple vehicles was perhaps first motivated by the work in distributed computing~\cite{Lynch96}, management science~\cite{DeGroot74}, and statistical physics~\cite{VicsekEtAl95}. In the control systems society, some pioneering works are generally referred to~\cite{Tsitsiklis84,TsitsiklisBertsekasAthans86}, where an asynchronous agreement problem was studied for distributed decision-making problems. Thereafter, some consensus algorithms were studied under various information-flow constraints~\cite{JadbabaieLinMorse03, FaxMurray04, SaberMurray04, RenBeard05_TAC,Moreau05}. There are several journal special issues on the related topics published after 2006, including the IEEE Transactions on Control Systems Technology (vol. 15, no. 4, 2007), Proceedings of the IEEE (vol. 94, no. 4, 2007), ASME Journal of Dynamic Systems, Measurement, and Control (vol. 129, no. 5, 2007), SIAM Journal of Control and Optimization (vol. 48, no.1, 2009), and International Journal of Robust and Nonlinear Control (vol. 21, no. 12, 2011). In addition, there are some recent reviews and progress reports given in the surveys~\cite{LeonardPLSFD07,RenBeardAtkins07_CSM,Murray07, SaberFaxMurray07_IEEE,ChebotarevAgaev09} and the books~\cite{Wu07_book,Shamma08,RenBeard08_Springer, BulloCortesMartinez09,Qu09,MesbahiEgerstedt10,RenCao11_Springer,BaiArcakWen11_Springer}, among others.

This article reviews some main results and recent progress in distributed multi-agent coordination, published in major control systems and robotics journals since 2006. For results before 2006, the readers are referred
to~\cite{LeonardPLSFD07,RenBeardAtkins07_CSM,Murray07,SaberFaxMurray07_IEEE}.

Specifically, this article reviews the recent research results in the following directions, which are not independent but actually may have overlapping to some extent:
\begin{itemize}
\item[{1.}]
Consensus and the like (synchronization, rendezvous). Consensus refers to the group behavior that all the agents asymptotically reach a certain common agreement through a local distributed protocol, with or without predefined common speed and orientation.
\item[{2.}]
Distributed formation and the like (flocking). Distributed formation refers to the group behavior that all the agents form a pre-designed geometrical configuration through local interactions with or without a common reference.
\item[{3.}]
Distributed optimization. This refers to algorithmic developments for the analysis and optimization of large-scale distributed systems.
\item[{4.}]
Distributed task assignment. This refers to the implementation of a task-assignment algorithm in a distributed fashion based on local information.
\item[{5.}]
Distributed estimation and control. This refers to distributed control design based on local estimation about the needed global information.
\end{itemize}

The rest of this article is organized as follows. In Section~\ref{sec:pre}, basic notations of graph theory and stochastic matrices are introduced. Sections~\ref{sec:consensus},~\ref{sec:formation},~\ref{sec:optimization},~\ref{sec:task}, and~\ref{sec:est} describe the recent research results and progress in consensus, formation control, optimization, task assignment, and estimation, respectively. Finally, the article is concluded by a short section of discussions with future perspectives.

\section{Preliminaries}\label{sec:pre}

This section introduces basic concepts and notations of graph theory and stochastic matrices.

\subsection{Graph Theory}\label{sec:graph}

For a system of $n$ connected agents, its network topology can be modeled as a directed graph denoted $\mathcal{G}=(\mathcal{V,W})$, where $\mathcal{V}=\{v_1,v_2, \cdots,v_n\}$ and $\mathcal{W} \subseteq\mathcal{V}\times\mathcal{V}$ are, respectively, the set of agents and the set of edges which directionally connect the agents together. Specifically, the directed edge denoted by an ordered pair $(v_i,v_j)$ means that agent $j$ can access the state information of agent $i$. Accordingly, agent $i$ is a neighbor of agent $j$. A directed path is a sequence of directed edges in the form of $(v_1,v_2),(v_2,v_3),\cdots\,$, with all $v_i\in \mathcal{V}$. A directed graph has a directed spanning tree if there exists at least one agent that has directed paths to
all other agents. The union of a set of directed graphs with the same set of agents, $\{\Gcal_{i_1},\cdots,\Gcal_{i_m}\}$, is a directed graph with the same set of agents and its set of edges is given by the union of the edge sets of all the directed graphs $\Gcal_{i_j}$, $j=1,\cdots,m$. A complete directed graph is a directed graph in which each pair of distinct agents is bidirectionally connected by an edge, thus there is a directed path from any agent to any other agent in the network.

Two matrices are frequently used to represent the network topology: the adjacency matrix $\Acal=[a_{ij}]\in\re^{n\times n}$ with $a_{ij}>0$ if $(v_j,v_i)\in \mathcal{W}$ and $a_{ij}=0$
otherwise, and the Laplacian matrix $\Lcal=[\ell_{ij}]\in \re^{n\times n}$ with $\ell_{ii}=\sum_{j=1}^na_{ij}$ and $\ell_{ij}=-a_{ij}$, $i\neq j$, which is generally asymmetric for directed graphs. The Laplacian $\Lcal$ has at least one single zero eigenvalue with a corresponding eigenvector $\onebf$ consisting of all numeric 1. Here and throughout, all matrices and vectors are assumed to have comparable dimensions unless otherwise indicated.

\subsection{Stochastic Matrices}\label{sec: stochas}

A nonnegative square matrix is called (row) stochastic matrix if its every row is summed up to one. The product of two stochastic
matrices is still a stochastic matrix. A row stochastic matrix $P\in\re^{n\times n}$ is called indecomposable and aperiodic if $\/\lim_{k\to\infty}P^k = \onebf y^T$ for some $y\in\re^{n}$~\cite{Wolfowitz63}.

\section{Consensus} \label{sec:consensus}

Consider a group of $n$ agents, each with single-integrator kinematics described by
\begin{equation}\label{eq:kinmatics}
  \dot{x}_i(t)=u_i(t),\quad i=1,\cdots,n,
\end{equation}
where $x_i(t)$ and $u_i(t)$ are, respectively, the state and the control input of the $\ith$ agent. A typical consensus control algorithm is designed as
\begin{equation}\label{eq:control}
  u_i(t)=\sum_{j=1}^n a_{ij}(t)[x_j(t)-x_i(t)],
\end{equation}
where $a_{ij}(t)$ is the $(i,j)$th entry of the corresponding adjacency matrix at time $t$. The main idea behind~\eqref{eq:control} is that each agent moves towards the weighted average of the states of its neighbors. Given the switching network pattern due to the continuous motions of the dynamic agents, coupling coefficients $a_{ij}(t)$ in~\eqref{eq:control}, hence the graph topologies, are generally time-varying. It is shown in~\cite{RenBeard05_TAC,Moreau05} that consensus is achieved if the underlying directed graph has a directed spanning tree in some jointly fashion in terms of a union of its time-varying graph topologies.

The idea behind consensus serves as a fundamental principle for the design of distributed multi-agent coordination algorithms. Therefore, investigating consensus has been a main research direction in the study of distributed multi-agent coordination. To bridge the gap between the study of consensus algorithms and many physical properties inherited in practical systems, it is necessary and meaningful to study consensus by considering many practical factors, such as actuation, control, communication, computation, and vehicle dynamics, which characterize some important features of practical systems. This is the main motivation to study consensus. In the following part of the section, an overview of the research progress in the study of consensus is given, regarding stochastic network topologies and dynamics, complex dynamical
systems, delay effects, the sampled-data framework,
asynchronous effects, quantization, convergence speed, and
finite-time convergence, mainly after 2006. Several milestone results prior to 2006 can be found in~\cite{DeGroot74,Tsitsiklis84,TsitsiklisBertsekasAthans86,JadbabaieLinMorse03, SaberMurray04,LinFrancisMaggiore05, RenBeard05_TAC,Moreau05}.

\subsection{Stochastic Network Topologies and Dynamics}~\label{sec:stochastic}
In multi-agent systems, the network topology among all vehicles plays a crucial role in determining consensus. The objective here is to explicitly identify necessary and/or sufficient conditions on the network topology such that consensus can be achieved under properly designed algorithms.


It is often reasonable to consider the case when the network topology is deterministic under ideal communication channels. Accordingly, main research on the consensus problem was conducted under a deterministic fixed/switching network topology. That is, the adjacency matrix $\Acal(t)$ is deterministic. Some other times, when considering random communication failures, random packet drops, and communication channel instabilities inherited in physical communication channels, it is necessary and important to study consensus problem in the stochastic setting where a network topology evolves according to some random distributions. That is, the adjacency matrix $\Acal(t)$ is stochastically evolving.

In the deterministic setting, consensus is said to be achieved if all agents eventually reach agreement on a common state. In the stochastic setting, consensus is said to be achieved \emph{almost surely} (respectively, \emph{in mean-square} or \emph{in probability}) if all agents reach agreement on a common state almost surely (respectively, in mean-square or in probability). Note that the problem studied in the stochastic setting is slightly different from that studied in the deterministic setting due to the different assumptions in terms of the network topology. Consensus over a stochastic network topology was perhaps first studied in~\cite{HatanoMesbahi05}, where some sufficient conditions on the network topology were given to guarantee consensus with probability one for systems with single-integrator kinematics~\eqref{eq:kinmatics}, where the rate of convergence was also studied. Further results for consensus under a stochastic network topology were reported in~\cite{Wu06,PorfiriStilwell07,AkarShorten08,SalehiJadbabaie08,HuangDNM10,HuangManton10,SalehiJadbabaie10_1,YinSunWang11,ZhangTian09}, where research effort was conducted for systems with single-integrator kinematics~\cite{Wu06,PorfiriStilwell07,AkarShorten08,SalehiJadbabaie08,HuangDNM10,HuangManton10,SalehiJadbabaie10_1,YinSunWang11} or double-integrator dynamics~\cite{ZhangTian09}. Consensus for single-integrator kinematics under stochastic network topology has been extensively studied in particular, where some general conditions for almost-surely consensus was derived~\cite{AkarShorten08,SalehiJadbabaie08,SalehiJadbabaie10_1}. Loosely speaking, almost-surely consensus for single-integrator kinematics can be achieved, \textit{i.e.,}~$x_i(t)-x_j(t)\to 0$ almost surely, if and only if the expectation of the network topology, namely, the network topology associated with expectation $E[\Acal(t)]$, has a directed spanning tree. It is worth noting that the conditions are analogous to that
in~\cite{Moreau05,RenBeard05_TAC}, but in the stochastic setting. In view of the special structure of the closed-loop systems concerning consensus for single-integrator kinematics, basic properties of the stochastic matrices play a
crucial role in the convergence analysis of the associated control algorithms. Consensus for double-integrator dynamics was studied in~\cite{ZhangTian09}, where the switching network topology is assumed to be driven by a Bernoulli process, and it was shown that consensus can be achieved if the union of all the graphs has a directed spanning tree. Apparently, the requirement on the network topology for double-integrator dynamics is a special case of that for single-integrator kinematics due to the difference nature of the final states (constant final states for single-integrator kinematics and possible dynamic final states for double-integrator dynamics) caused by the substantial dynamical difference. It is still an open question as if some general conditions (corresponding to some specific algorithms) can be found for consensus with double-integrator dynamics.

In addition to analyzing the conditions on the network topology such that consensus can be achieved, a special type of consensus algorithm, the so-called gossip algorithm~\cite{BoydGPS06,DimakisKMRS10}, has been used to achieve consensus in the stochastic setting. The gossip algorithm can always guarantee consensus almost surely if the available pairwise communication channels satisfy certain conditions (such as a connected graph or a graph with a directed spanning tree). The way of network topology switching does not play any role in the consideration of consensus.

The current study on consensus over stochastic network topologies has shown some interesting results regarding: (1) consensus algorithm design for various multi-agent systems, (2) conditions of the network topologies on consensus, and (3) effects of the stochastic network topologies on the convergence rate. Future research on this topic includes, but not limited to, the following two directions: (1) when the network topology itself is stochastic, how to determine the probability of reaching consensus almost surely? (2) compared with the deterministic network topology, what are the advantages and disadvantages of the stochastic network topology, regarding such as robustness and convergence rate?

As is well known, disturbances and uncertainties often exist in networked systems, for example, channel noise, communication noise, uncertainties in network parameters, etc. In addition to the stochastic network topologies discussed above, the effect of stochastic disturbances~\cite{LiZhang09,LiDuanChen11} and uncertainties~\cite{KimShimSeo11} on the consensus problem also needs investigation. Study has been mainly devoted to analyzing the performance of consensus algorithms subject to disturbances and to presenting conditions on the uncertainties such that consensus can be achieved. In addition, another interesting direction in dealing with disturbances and uncertainties is to design distributed local filtering algorithms so as to save energy and improve computational efficiency. Distributed local filtering algorithms play an important role and are more effective than traditional centralized filtering algorithms for multi-agent systems. For example, in~\cite{SaberCDC05a,SaberCDC05b,SaberCDC07} some distributed Kalman filters are designed to implement data fusion. In~\cite{YuTSMCB09}, by analyzing consensus and pinning control in synchronization of complex networks, distributed consensus filtering in sensor networks is addressed. Recently, Kalman filtering over a packet-dropping network is designed through a probabilistic approach \cite{ShiTAC10}. Today, it remains a challenging problem to incorporate both dynamics of consensus and probabilistic (Kalman) filtering into a unified framework.


\subsection{Complex Dynamical Systems}
Since consensus is concerned with the behavior of a group of vehicles, it is natural to consider the system dynamics for practical vehicles in the study of the consensus problem. Although the study of consensus under various system dynamics is due to the existence of complex dynamics in practical systems, it is also interesting to observe that system dynamics play an important role in determining the final consensus state. For instance, the well-studied consensus of multi-agent systems with single-integrator kinematics often converges to a constant final value instead. However, consensus for double-integrator dynamics might admit a dynamic final value (i.e., a time function). These important issues motivate the study of consensus under various system dynamics.

As a direct extension of the study of the consensus problem for systems with simple dynamics, for example, with single-integrator kinematics or double-integrator dynamics, consensus with general linear dynamics was also studied recently~\cite{QuWangHull08,Tuna09,ScardoviSepulchre09,SeoShimBack09,LiDCH10,YangRWS11}, where research is mainly devoted to finding feedback control laws such that consensus (in terms of the output states) can be achieved for general linear systems
\begin{equation}~\label{eq:general-linear-dyn}
\dot{x}_i=Ax_i+Bu_i,\quad y_i=Cx_i,
\end{equation}
where $A,~B,~\text{and}~C$ are constant matrices with compatible sizes. Apparently, the well-studied single-integrator kinematics and double-integrator dynamics are special cases of~\eqref{eq:general-linear-dyn} for properly choosing $A,~B,~\text{and}~C$.


As a further extension, consensus for complex systems has also been extensively studied. Here, the term {\it consensus for complex systems\/} is used for the study of consensus problem when the system dynamics are
nonlinear~\cite{StanSepulchre07,LinFrancisMaggiore07, DimarogonasKyriakopoulos08,ZhouLuLu06,ZhouLuLu08,ChopraSpong09, YuCCLHW09,ChenCXLY09,HouChengTan09,DelellisDibernandoGarofalo09, ShiHong09,SarletteSepulchreLeonard09,BaiArcakWen08,ChungSlotine09, NairLeonard08,YaoGuanHill09,ZhaoHillLiu09,MenonEdwards09,WangWXG10,LuHoCao10,EstradaGagoCaporossi10,DasLewis10,SuCWL11,YuChenCao11,
ScardoviArcakSontag10,LiJiang09,AbdessameudTayebi09,ChenLewis11,ZhangMHW10,LuHo10,WangWang11,PowelMorgansen11,WangXiaoWang10,Ren09_IJC} or with nonlinear consensus algorithms~\cite{HuiHaddadBhat08,HuiHaddad08,Cortes06}. Examples of the nonlinear system dynamics studied in the consensus problem include:
\begin{itemize}
\item Nonlinear oscillators~\cite{ChopraSpong09}. The dynamics are often assumed to be governed by the Kuramoto equation
\begin{equation}
\dot{\theta}_i=\omega_i+\frac{K}{N}\sum_{j=1}^N \sin(\theta_j-\theta_i),\quad i=1,2,\cdots,N,
\end{equation}
where $\theta_i$ and $\omega_i$ are, respectively, the phase and natural frequency of the $i$th oscillator, $N$ is the number of oscillators, and $K$ is the control gain. Generally, the control gain $K$ plays a crucial role in determining the synchronizability of the network.
\item Complex
networks~\cite{ZhouLuLu06,YuChenLu09,WangWang11,YaoGuanHill09,ZhaoHillLiu09,MenonEdwards09,WangWXG10,LuHoCao10,ZhangMHW10,LuHo10,WangXiaoWang10}. The dynamics are typically represented as
\begin{align}
\dot{x}_{i}(t)=& f(x_{i}(t))+c \sum\limits_{j=1,j\neq
i}^{N}a_{ij}(t)\Gamma(x_{j}(t)-x_{i}(t)),\notag
\\&~~~~~~~~~~~~~i=1,2,\cdots,N,\label{synchronization}
\end{align}
where $x_i=\left(x_{i1},x_{i2},\cdots,x_{in}\right)^T\in \mathbb
R^n$ is the state vector of the $i$th node,
$f\,:\,\mathbb{R}^{n}\mapsto
\mathbb{R}^{n}$ is a nonlinear vector function, $c$ is the
overall coupling strength, $A(t)=\left[a_{ij}(t)\right]$ is the outer coupling matrix with $a_{ij}(t)=1$ if node $i$ and node $j$ are connected at time $t$ but otherwise $a_{ij}(t)=0$, with $a_{ii}(t)=k_i$ (degree of node $i$), and $\Gamma$ is
a general inner coupling matrix describing the inner interactions between different state components of agents. It is easy to see that model (\ref{eq:kinmatics}) with control input (\ref{eq:control}) is a special case of (\ref{synchronization}) with $f=0$.
\item Nonholonomic mobile
robots~\cite{DimarogonasKyriakopoulos07,DimarogonasKyriakopoulos08,LiJiang09,PowelMorgansen11}. The dynamics are described by
\begin{align}\label{eq:nonholonomic}
\dot{x}_i=u_i\cos\theta_i, \dot{y}_i=u_i\sin\theta_i, \dot{\theta}_i=\omega_i, i=1,\cdots,N,
\end{align}
where $[x_i,y_i]$ denotes the location of the $i$th agent, and $u_i$ and $\omega_i$ denote, respectively, its translational and rotational velocity. Note that there are three states and two control inputs. Therefore, the dynamics for nonholonomic mobile
robots are underactuated. This poses substantial difficulties in designing proper consensus algorithms with corresponding stability analysis.
\item Rigid bodies and the like~\cite{NairLeonard08, SarletteSepulchreLeonard09,BaiArcakWen08,ChungSlotine09,ChenLewis11,Ren09_IJC,AbdessameudTayebi09}. One typical (but not unique) description of the dynamics is
    \begin{equation}~\label{eq:rigid-body}
    M_i(q_i)\ddot{q}_i+C_i(q_i,\dot{q}_i)\dot{q}_i+g_i(q_i)=\tau_i,~ i=1,\cdots,N,
    \end{equation}
    where $q_i\in \re^p$ is the vector of generalized coordinates, $M_i(q_i)\in\re^{p\times p}$ is the symmetric
positive-definite inertia matrix, $C_i(q_i,\dot{q}_i)\dot{q}_i\in\re^p$ is the vector of Coriolis and
centrifugal torques, $g_i(q_i)$ is the vector of gravitational torques, and $\tau_i\in\re^p$ is the
vector of torques produced by the actuators associated with the $i$th agent. In practice, the dynamics of many mechanical systems are similar to~\eqref{eq:rigid-body}. A notable property regarding the dynamics of rigid bodies is that $\dot{M}_i(q_i)-2C_i(q_i,\dot{q}_i)$ is skew-symmetric (\textit{i.e.,} $z^T[\dot{M}_i(q_i)-2C_i(q_i,\dot{q}_i)]z=0$ for all $z\in\re^p$), which plays a crucial role in finding Lyapunov functions and the subsequent stability analysis.
\end{itemize}
Although the aforementioned system dynamics are different from the well-studied single-integrator kinematics and double-integrator dynamics, the main
research problem is same, namely, to drive all agents to some common states through local interactions among agents. Similarly to the consensus algorithms proposed for systems with simple dynamics, the consensus algorithms used for these complex models are also based on a weighted average of the state differences, with some additional terms if necessary.  Main research work has been conducted to design proper control algorithms and derive necessary and/or sufficient conditions such that consensus can be achieved ultimately.


Note that although the objective is same, \emph{i.e.,} to guarantee reaching agreement on some final states, the problem is more complicated due to the nonlinearity of the closed-loop systems. In addition, most properties of stochastic matrices cannot be directly applied to their convergence analysis. The main techniques used in their stability analysis include dissipativity theory~\cite{StanSepulchre07}, nonsmooth
analysis~\cite{DimarogonasKyriakopoulos07,HuiHaddad08,Cortes06}, and especially Lyapunov functions~\cite{DimarogonasKyriakopoulos07, ChopraSpong09,NairLeonard08,SarletteSepulchreLeonard09,BaiArcakWen09, ChungSlotine09}.

One particular interesting topic is synchronization in complex networks which has been widely investigated in the past decade \cite{WangChenTCASI02,WangChenIJBC02}.
%
Mathematically, the definitions for synchronization in complex networks and consensus in multi-agent systems are very similar, so to differentiate these two definitions and promote research exchanges in these two topics, their differences are briefly summarized below.

1) {\it Different Asymptotic States} (Nonlinear Dynamics versus Linear Dynamics). In the studies of synchronization in complex networks, researchers focus on synchronization with self-nonlinear dynamics where each single system is unstable and thus the final asymptotic synchronization state is typically time-varying \cite{WangChenTCASI02,YuSIAM08}. However, in the investigations of multi-agent systems, the individual self-dynamics on each system is usually linear or zero and therefore the asymptotic consensus state is usually a constant \cite{JadbabaieLinMorse03,SaberMurray04}.

2) {\it Different Focuses} (Known Connectivity versus Time-varying Distributed Protocol). In synchronization of complex networks, the aim is to reveal how the network structure, which is known in priori, affects the nonlinear collective dynamics \cite{WangChenTCASI02,WangChenIJBC02}, while the aim of consensus in multi-agent systems is to figure out how the designed distributed local protocol concerning mobile time-varying network structure affects the consensus behavior \cite{JadbabaieLinMorse03,SaberMurray04,RenBeard05_TAC}.

3) {\it Different Approaches} (Lyapunov Method versus Stochastic Matrix Theory). Since both complex networks and multi-agent systems are networked systems, algebraic graph theory \cite{Godsil01} is a common approach to use. Because of the nonlinear terms in synchronization of complex networks, Lyapunov function method is usually used together with matrix theory \cite{ZhouLuLu06,YuChenLu09,YuSIAM08}. In order to show consensus in multi-agent systems with time-varying network structures, stochastic matrix theory \cite{Tsitsiklis84,TsitsiklisBertsekasAthans86,JadbabaieLinMorse03,RenBeard05_TAC} and convexity analysis \cite{Moreau05} are often applied.

4) {\it Different Inner Matrices $\Gamma$} (General Inner Matrix versus Particular Inner Matrix). In the typical simple consensus model, the inner matrices $\Gamma$ are usually an identity matrix and a rank-one matrix $\left(
                                                 \begin{array}{cc}
                                                   0 & 1 \\
                                                   0 & 0 \\
                                                 \end{array}
                                               \right)$
for multi-agent systems with single-integrator kinematics \cite{SaberMurray04} and  double-integrator dynamics \cite{RenAtkins07,YuChenCao10,YuTSMCB10}, respectively. In consensus models with higher-order dynamics \cite{YuTCASI11}, the inner matrix is similar. However, the inner matrix in system (\ref{synchronization}) is a general one.

In summary, synchronization in complex networks focuses on nonlinear dynamics while consensus in multi-agent systems focuses on distributed cooperative control, and thus different approaches are utilized.

The current research on consensus with complex systems focuses on fully-actuated systems although consensus for nonholonomic mobile robots~\cite{DimarogonasKyriakopoulos07}, which are typical underactuated systems. Note that many mechanical devices are described by systems with underactuation. Therefore, it is important to develop appropriate consensus algorithms for underactuated systems.

\subsection{Delay Effects}
Time delay appears in almost all practical systems due to several reasons: (1) limited communication speed when information transmission exists; (2) extra time required by the sensor to get the measurement information; (3) computation time required for generating the control inputs; and (4) execution time required for the inputs being acted. In general, time delay reflects an important property inherited in practical systems due to actuation, control, communication, and computation.

Knowing that time delay might degrade the system performance or even destroy the system stability, studies have been conducted to investigate its effect on system performance and stability. A well-studied consensus algorithm for~\eqref{eq:kinmatics} is given in~\eqref{eq:control}, where it is now assumed that time delay exists. Two types of time delays, \emph{communication delay} and \emph{input delay}, have been considered in the literature. Communication delay accounts for the time for transmitting information from origin to destination. More precisely, if it takes time $T_{ij}$ for agent $i$ to receive information from agent $j$, the closed-loop system of~\eqref{eq:kinmatics} using~\eqref{eq:control} under a fixed network topology becomes
\begin{equation}\label{eq:comm}
    \dot{x}_i(t)=\sum_{j=1}^n a_{ij}(t)[x_j(t-T_{ij})-x_i(t)].
\end{equation}
An interpretation of~\eqref{eq:comm} is that at time $t$, agent $i$ receives information from agent $j$ and uses data $x_j(t-T_{ij})$ instead of $x_j(t)$ due to the time delay. Note that agent $i$ can get its own information instantly, therefore, input delay can be considered as the summation of computation time and execution time. More precisely, if the input delay for agent $i$ is given by $T^p_i$, then the closed-loop system of~\eqref{eq:kinmatics} using~\eqref{eq:control} becomes
\begin{equation}\label{eq:input}
    \dot{x}_i(t)=\sum_{j=1}^n a_{ij}(t)[x_j(t-T_i^{p})-x_i(t-T_i^{p})].
\end{equation}
Clearly,~\eqref{eq:comm} refers to the case when only communication delay is considered while~\eqref{eq:input} refers to the case when only input delay is considered. It should be emphasized that both communication delay and input delay might be time-varying and they might co-exist at the same time.

In addition to time delay, it is also important to consider packet drops in exchanging state information. Fortunately, consensus with packet drops can be considered as a special case of consensus with time delay, because re-sending packets after they were dropped can be easily done but just having time delay in the data transmission channels.

Thus, the main problem involved in consensus with time delay is to study the effects of time delay on the convergence and performance of consensus, referred to as {\em
consensusability}~\cite{MaZhang10}.


Because time delay might affect the system stability, it is important to study under what conditions consensus can still be guaranteed even if time delay exists.
In other words, can one find conditions on the time delay such that consensus can be achieved? For this purpose, the effect of time delay on the consensusability of~\eqref{eq:kinmatics} using~\eqref{eq:control} was investigated. When there exists only (constant) input delay, a sufficient condition on the time delay to guarantee consensus under a fixed undirected interaction graph is presented in~\cite{SaberMurray04}. Specifically, an upper bound for the time delay is derived under which consensus can be achieved. This is a well-expected result because time delay normally degrades the system performance gradually but will not destroy the system stability unless the time delay is above a certain threshold. Further studies can be found in, e.g.,~\cite{WangSlotine06_1,XiaoWang08,XiaoWang08_TAC,SunWangXie,LinJiaLi08,TianLiu08,SunWang09,MunzPapachristodoulouAllgower10,LiuLuChen10,PapachristodoulouJadbabaieMunz10,LinJia11}, which demonstrate that for~\eqref{eq:kinmatics} using~\eqref{eq:control}, the communication delay
does not affect the consensusability but the input delay does.
In a similar manner, consensus with time delay was studied for systems with different dynamics, where the dynamics~\eqref{eq:kinmatics} are replaced by other more complex ones, such as double-integrator dynamics~\cite{TianLiu09,LinJia09,YangFang10,YuChenCao10,ZhangTian10,LinJia10,XiaoWangChen10,MunzPapachristodoulouAllgower11,QinGaoZheng11}, complex networks~\cite{CaoChenLi08,LiangWLL08,YaoWGX09,WangZWY10}, rigid bodies~\cite{ChopraSpongLozano08,NunoOBH11}, and general nonlinear dynamics~\cite{BokharaieMasonVerwoerd10}.



In summary, the existing study of consensus with time delay mainly focuses on analyzing the stability of consensus algorithms with time delay for various types of system dynamics, including linear and nonlinear dynamics. Generally speaking, consensus with time delay for systems with nonlinear dynamics is more challenging. For most consensus algorithms with time delays, the main research question is to determine an upper bound of the time delay under which time delay does not affect the consensusability. For communication delay, it is possible to achieve consensus under a relatively large time delay threshold. A notable phenomenon in this case is that the final consensus state is constant. Considering both linear and nonlinear system dynamics in consensus, the main tools for stability analysis of the closed-loop systems include matrix theory~\cite{XiaoWang08,XiaoWang08_TAC},
Lyapunov functions~\cite{LiangWLL08,CaoChenLi08}, frequency-domain approach~\cite{TianLiu08}, passivity~\cite{YaoWGX09}, and the contraction principle~\cite{WangSlotine06_1}.

Although consensus with time delay has been studied extensively, it is often assumed that time delay is either constant or random. However, time delay itself might obey its own dynamics, which possibly depend on the communication distance, total computation load and computation capability, etc. Therefore, it is more suitable to represent the time delay as another system variable to be considered in the study of the consensus problem. In addition, it is also important to consider time delay and other physical constraints simultaneously in the study of the consensus problem.

\subsection{Sampled-data Framework}
The previous three subsections describe the main research work in the study of the consensus problem. The following introduces a few other aspects, namely, sampled-data framework, quantization, asynchronous effect, convergence speed, and finite-time convergence, that have been considered in the consensus problem as well. Among these topics, sampled-data framework, quantization, and asynchronous effects are considered due to some physical limitations in practical systems while convergence speed and finite-time convergence are concerned with the performance for some proposed consensus algorithms.

Due to the limitations in the measurement and control units, it is often impossible to acquire information measurements at an arbitrarily fast
speed and to execute the control inputs instantaneously. Accordingly, the closed-loop systems are modeled in a hybrid fashion. That is, the system plant is described in a continuous-time setting while the measurements and control inputs are described in a piecewise constant fashion. For instance, in a sampled-data setting,~\eqref{eq:control} becomes
\begin{align}\label{eq:discrete-time}
u_i(t)=u_i(kT)=&\sum_{j=1}^n a_{ij}(kT)[x_j(kT)-x_i(kT)]
\end{align}
for $kT\leq t<(k+1)T$, where $T$ is the sampling period and $k$ is the discrete-time index. Essentially,~\eqref{eq:discrete-time} is a zero-order-hold version of $u_i(t)$ in the sense that the control inputs remain unchanged during each sampling period.
Under this circumstance, consensus is studied in a sampled-data framework, called {\it sampled-data consensus}, which reflects the limitations inherited in physical measurement and control units. Meanwhile, it is also important to point out that the sampled-data consensus algorithms require much less information exchange and computational power than the continuous-time consensus algorithms. Accordingly, consensus under the sampled-data framework deserves certain consideration.


Sampled-data consensus was investigated in, e.g.,~\cite{HayakawaMatsuzawaHara06,XieLWJ09_1,XieLWJ09_2,GaoWXW09, GaoWang10,CaoRen09_IJRNC,CaoRen09_IJC,YuZCRC11,GaoWang11,LiuXieWang10,ZhangTian10}. Consensus for systems with single-integrator kinematics~\eqref{eq:kinmatics} was studied under a sampled-data framework with a fixed or a switching network topology, in~\cite{XieLWJ09_1,XieLWJ09_2}, where some necessary and/or sufficient conditions were presented to guarantee achieving consensus. Sampled-data consensus of systems with double-integrator kinematics was studied under fixed or switching network topologies in~\cite{HayakawaMatsuzawaHara06, GaoWXW09,GaoWang10,CaoRen09_IJRNC,CaoRen09_IJC,ZhangTian10,YuZCRC11,GaoWang11,LiuXieWang10}. Due to the fact that an infinitely large sampling period will cause no information exchange among the agents, the main research question is to find conditions on the sampling period $T$, which might be time-varying, such that consensus can be achieved. The conditions on the network topology for the sampled-data closed-loop systems are mostly similar to that for the continuous-time closed-loop systems.
Note that the existing research on consensus in a sampled-data framework mainly focuses on the simple system dynamics and thus the closed-loop system can be represented in terms of a linear matrix equation. The corresponding network stability can be analyzed by investigating the properties of the system matrices constructed based on the proposed consensus algorithms and the given network topology. Various approaches, including Lyapunov functions~\cite{HayakawaMatsuzawaHara06,GaoWang10}, matrix
theory~\cite{CaoRen09_IJRNC,ZhangTian10,GaoWang11,LiuXieWang10}, stochastic matrices~\cite{CaoRen09_IJC}, and
linear matrix inequalities~\cite{GaoWXW09,GaoWang10}, have been adopted, and some necessary and/or sufficient conditions have been derived for guaranteeing sampled-data consensus.

It is natural to consider the sampled-data effect for consensus with general linear or nonlinear dynamics (see,~\textit{e.g.,}~\cite{WenDYC11}). In addition, it is meaningful to consider the case when all vehicles do not necessarily share the same sampling period or the sampling period is not necessarily constant. Accordingly, it is expected that a careful design of the sampling periods (associated with the proposed algorithms) might lead to the optimization of the closed-loop systems under the proposed algorithms subject to certain cost functions, such as maximum convergence rate and minimum total information exchange. In another word, it is intriguing to move from analysis to design when investigating the consensus problem in a sampled-data framework.


\subsection{Asynchronous Effects}

In most existing research of the consensus problem, it is assumed that all agents update their states synchronously, which requires a synchronized clock for the whole group of agents. However, such a synchronized clock might not exist in real applications. This motivates the design of consensus algorithms in an asynchronous fashion; that is, each agent updates its own states regardless of the update times of other agents.

In most studies of the asynchronous consensus problem for networked systems, due to the intrinsic technical difficulties, usually only single-integrator kinematics~\eqref{eq:kinmatics} and double-integrator dynamics are considered. In~\cite{XiaoWang08_TAC}, such an asynchronous consensus problem with time delay was investigated by utilizing some basic properties of stochastic matrices. Similarly
in~\cite{CaoMorseAnderson08}, the asynchronous consensus problem was studied by using matrix theory and graph theory, and in~\cite{FangAntsaklis08}, by employing the paracontracting theory. In~\cite{GaoWang10_IJC}, the authors studied the asynchronous consensus problem for double-integrator dynamics and presented sufficient conditions to guarantee consensus, where a condition based on linear matrix inequalities was given.

Note that consensus in an asynchronous fashion has been considered mainly for single-integrator kinematics and double-integrator dynamics but not for other system dynamics. For certain linear systems, it might be expected that asynchronous communication does not affect the consensusability as shown in~\cite{XiaoWang08_TAC,CaoMorseAnderson08} for single-integrator kinematics. However, a similar conclusion may not hold for systems with general dynamics, especially nonlinear dynamics. It is important to quantify the effects of the asynchronous communication on the consensus problem.

\subsection{Quantization}

Quantized consensus has been studied recently with motivation from digital signal processing. Here, quantized consensus refers to consensus when the measurements are digital rather than analog therefore the information received by each agent is not continuous and might have been truncated due to digital finite precision constraints. Roughly speaking, for an analog signal $s$, a typical quantizer with an accuracy parameter $\delta$, also referred to as quantization step size, is described by
$
Q(s) = q(s,\delta),
$
where $Q(s)$ is the quantized signal and $q(\cdot,\cdot)$ is the associated quantization function. For instance~\cite{KashyapBasarSrikant07}, a quantizer rounding a signal $s$ to its nearest integer can be expressed as
$
Q(s) = n,~\text{if}~s\in [(n-1/2)\delta,(n+1/2)\delta],~ n\in \Zcal,
$
where $\Zcal$ denotes the integer set. Note that the types of quantizers might be different for different systems, hence $Q(s)$ may differ for different systems.
Due to the truncation of the signals received, consensus is now considered achieved if the maximal state difference is not larger than the accuracy level associated with the whole system. A notable feature for consensus with quantization is that the time to reach consensus is usually finite. That is, it often takes a finite period of time for all agents' states to converge to an accuracy interval. Accordingly, the main research is to investigate the convergence time associated with the proposed consensus algorithm.

Quantized consensus was probably first studied in~\cite{KashyapBasarSrikant07}, where a quantized gossip algorithm was proposed and its convergence was analyzed. In particular, the bound of the convergence time for a complete graph was shown to be polynomial in the network size. In~\cite{CarliBullo09}, coding/decoding strategies were introduced to the quantized consensus algorithms, where it was shown that the convergence rate depends on the accuracy of the quantization but not the coding/decoding schemes. In~\cite{LavaeiMurray09_1, LavaeiMurray09_2}, quantized consensus was studied via the gossip algorithm, with both lower and upper bounds of the convergence time derived in terms of the network size. Further results regarding quantized consensus were reported in~\cite{YuanXZC10,NedicOOT09,FranceschelliGiuaSeatzu10,LiFXZ11,ZhuMartinez11,CarliCFG11,DimarogonasJohansson10,CarliFFZ10}, where the main research was also on the convergence time for various proposed quantized consensus algorithms as well as the quantization effects on the convergence time. It is intuitively reasonable that the convergence time depends on both the quantization level and the network topology. It is then natural to ask if and how the quantization methods affect the convergence time. This is an important measure of the robustness of a quantized consensus algorithm (with respect to the quantization method).

Note that it is interesting but also more challenging to study consensus for general linear/nonlinear systems with quantization. Because the difference between the truncated signal and the original signal is bounded, consensus with quantization can be considered as a special case of one without quantization when there exist bounded disturbances. Therefore, if consensus can be achieved for a group of vehicles in the absence of quantization, it might be intuitively correct to say that the differences among the states of all vehicles will be bounded if the quantization precision is small enough. However, it is still an open question to rigorously describe the quantization effects on consensus with general linear/nonlinear systems.

\subsection{Convergence Speed}
In addition to the study on the consensus problem with physical constraints mentioned in the previous subsections, it is also important to study the control performance of the consensus problem. From the control's perspective, it is natural to propose proper control algorithms and analyze the stability, and to optimize the proposed control algorithms under certain control performance indexes. In this subsection, the convergence speed problem is reviewed, which is an important performance measure for consensus algorithms.

For~\eqref{eq:kinmatics} using algorithm~\eqref{eq:control} in a connected undirected graph, the worst-case convergence speed was shown in~\cite{SaberMurray04} to be the Laplacian spectral gap
$
\min_{ X\neq\zerobf,\onebf^TX=0}\frac{X^T \Lcal X}{\norm{X}^2}=\lambda_2,
$
where $\zerobf$ is an all-zero column vector, $X=[x_1,\cdots, x_n]^T$, $\Lcal$ is the Laplacian matrix with $\lambda_2$ being smallest nonzero eigenvalue. Here, one should recall that the smallest eigenvalue of a Laplacian matrix for a connected undirected graph is zero and all the other eigenvalues are positive.

In order to increase the convergence speed, the above spectral gap should be enlarged. For this purpose, an iterative algorithm was proposed in~\cite{KimMesbahi06} to maximize the above spectral gap, by employing a semidefinite programming solver.

Other than the smallest nonzero eigenvalue of the Laplacian matrix, another commonly used measure for the convergence speed is the following ratio, introduced in~\cite{XiaoBoyd04,OlshevskyTsitsiklis09}:
\begin{equation}\label{eq:convspeed}
  \rho=\lim_{t\to\infty,X(t)\neq
  X^\star}\,\left(\frac{\norm{X(t)-X^\star}}
  {\norm{X(0)-X^\star}}\right)^{1/t},
\end{equation}
where $X^\star$ represents the final equilibrium given by $\sigma \onebf$, where $\sigma$ is a constant.

In~\cite{XiaoBoyd04}, this problem of finding the fastest convergence speed was casted into a semidefinite programming problem. Furthermore, the convergence speed defined by~\eqref{eq:convspeed} was studied in both deterministic and stochastic settings. In the
deterministic setting, it was studied in~\cite{AngeliBliman09, OlshevskyTsitsiklis09,AngeliBliman08} with estimation of lower bounds. In the stochastic setting, this problem was studied in~\cite{HatanoMesbahi05,AkarShorten08,ZhouWang09}, with a per-step convergence factor introduced and discussed in~\cite{ZhouWang09}, which itself can be considered a measure of the convergence speed.


The existing study mainly focuses on the analysis of the convergence speed under various network topologies and optimization of the convergence speed for certain given network topologies. Considering the fact that consensus under different network topologies may demonstrate different convergence speeds, a natural question is how to design an optimal (switching) network topology with proper adjacency matrix such that optimal convergence speed can be achieved.

\subsection{Finite-time Convergence}
As an extension of the study of convergence speed for the consensus problem, \emph{finite-time consensus}, reaching consensus in a finite time, has also been studied recently. Compared with most existing research on the consensus problem, finite-time consensus demonstrates a disturbance rejection property and robustness against uncertainties. In addition, due to the finite-time convergence, it is often possible to decouple the consensus problem from other control objectives when they are considered simultaneously.


For a group of $n$ agents, e.g. with dynamics~\eqref{eq:kinmatics}, the objective is to design
$u_i(t)$ such that $x_i(t)=x_j(t)$ for $t\geq\overline{T}$, where $\overline{T}$ is a constant. Here, $\overline{T}$ is called the consensus time.

Finite-time consensus for networked systems with single-integrator kinematics~\eqref{eq:kinmatics} in the continuous-time setting was solved in~\cite{HuiHaddadBhat08,Cortes06,XiaoWCG09,JiangWang09,CaoRenMeng10,WangXiao10,Hui11}. Finite-time consensus for networked systems with double-integrator dynamics in the continuous-time setting was studied in~\cite{WangHong08}. An important common characteristic of the proposed algorithms is the use of the signum function. It is well known that linear consensus algorithms can normally guarantee asymptotic convergence, but not finite-time convergence. On the contrary, the finite-time consensus algorithms developed in~\cite{HuiHaddadBhat08,Cortes06,WangHong08,XiaoWCG09,JiangWang09,CaoRenMeng10,Hui11}, which utilize the signum function, are able to do so.

Note that the existing research on finite-time consensus mainly focuses on systems with simple dynamics, such as single-integrator kinematics and double-integrator dynamics, in the continuous-time setting. Because many practical systems are better and more proper to be described by general linear/nonlinear dynamics, it is natural to study finite-time consensus for systems with general linear/nonlinear dynamics in the future. In addition, it is interesting to study finite-time consensus in the discrete-time setting. Some recent work on this topic can be found in~\cite{SundaramHadjicostis07,YuanSBSG11}, where the objective is to compute the final consensus value for all agents in a finite number of steps.

\subsection{Remarks}
In summary, the existing research on the consensus problem has covered a number of physical properties for practical systems and control performance analysis. However, the study of the consensus problem covering multiple physical properties and/or control performance analysis has been largely ignored. In other words, two or more problems discussed in the above subsections might need to be taken into consideration simultaneously when studying the consensus problem. In addition, consensus algorithms normally guarantee the agreement of a team of agents on some common states without taking group formation into consideration. To reflect many practical applications where a group of agents are normally required to form some preferred geometric structure, it is desirable to consider a task-oriented formation control problem for a group of mobile agents, which motivates the study of formation control presented in the next section.

\section{Formation Control}\label{sec:formation}

Compared with the consensus problem where the final states of all agents typically reach a singleton, the final states of all agents can be more diversified under the formation control scenario. Indeed, formation control is more desirable in many practical applications such as formation flying, cooperative transportation, sensor networks, as well as combat intelligence, surveillance, and reconnaissance. In addition, the performance of a team of agents working cooperatively often exceeds the simple integration of the performances of all individual agents. For its broad applications and advantages, formation control has been a very active research subject in the control systems community, where a certain geometric pattern is aimed to form with or without a group reference. More precisely, the main objective of formation control is to coordinate a group of agents such that they can achieve some desired formation so that some tasks can be finished by the collaboration of the agents. Generally speaking, formation control can be categorized according to the group reference. Formation control without a group reference, called {\it formation producing}, refers to the algorithm design for a group of agents to reach some pre-desired geometric pattern in the absence of a group reference, which can also be considered as the control objective. Formation control with a group reference, called {\it formation tracking}, refers to the same task but following the predesignated group reference. Due to the existence of the group reference, formation tracking is usually much more challenging than formation producing and control algorithms for the latter might not be useful for the former. As of today, there are still many open questions in solving the formation tracking problem.

The following part of the section reviews and discusses recent research results and progress in formation control, including formation producing, formation tracking, and connectivity maintenance for consensus and formation control, mainly accomplished after 2006. Several milestone results prior to 2006 can be found in~\cite{Reynolds87,LeonardFiorelli01,OgrenFiorelliLeonard04}.

\subsection{Formation Producing}\label{sec:formatioProducing}

The existing work in formation producing aims at analyzing the formation behavior under certain control laws, along with stability analysis.

\subsubsection{Matrix Theory Approach}

Due to the nature of multi-agent systems, matrix theory has been frequently used in the stability analysis of their distributed coordination.

Note that consensus input to each agent (see e.g.,~\eqref{eq:control}) is essentially a weighted average of the differences between the states of the agent's neighbors and its own. As an extension of the consensus algorithms, some coupling matrices were introduced here to offset the corresponding control inputs by some angles~\cite{PavoneFrazzoli07,Ren09_TAC,LinJia10_SCL,LinQLR11}. For example, given~\eqref{eq:kinmatics}, the control input~\eqref{eq:control} is revised as
\begin{equation}\label{eq:control-rotation}
  u_i(t)=\sum_{j=1}^n a_{ij}(t)C[x_j(t)-x_i(t)],
\end{equation}
where $C$ is a coupling matrix with compatible size. If $x_i\in\re^3$, then $C$ can be viewed as the 3-D rotational matrix. The main idea behind~\eqref{eq:control-rotation} is that the original control input for reaching consensus is now rotated by some angles. The closed-loop system can be expressed in a vector form, whose stability can be determined by studying the distribution of the eigenvalues of a certain transfer matrix. Main research work was conducted in~\cite{PavoneFrazzoli07,Ren09_TAC,LinJia10_SCL,LinQLR11} to design proper algorithms and analyze the collective motions for systems with single-integrator kinematics and double-integrator dynamics. Analogously, the collective motions for a team of
nonlinear self-propelling agents were shown to be affected by the coupling strength among the agents, the time delay, the
noise, and the initial states~\cite{RomeroForgostonSchwartz12}.


Note that the collective motions for nonholonomic mobile robots were also studied recently in, e.g.,~\cite{SepulchrePaleyLeonard07,SepulchrePaleyLeonard08}. Although the study in~\cite{PavoneFrazzoli07,Ren09_TAC,LinJia10_SCL,LinQLR11} is different from that in~\cite{SepulchrePaleyLeonard07,SepulchrePaleyLeonard08}, similarities exist in the sense that all agents will not move to the weighted average of the states of neighboring agents, but to some offsetted state. Noticeably, the offsetted state in~\cite{PavoneFrazzoli07,Ren09_TAC,LinJia10_SCL,LinQLR11} is properly designed, yet the one in~\cite{SepulchrePaleyLeonard07,SepulchrePaleyLeonard08} is induced by some special nonlinear system dynamics.

In the study of formation producing with linear closed-loop systems, it is observed that the associated system matrix has two interesting properties: (1) the existence of at least one zero eigenvalue, and (2) the existence of at least one pair of eigenvalues on the imaginary axis. The two properties play an important role in the formation producing problem under a fixed network topology. However, the two properties might not be able to solve the formation producing problem under a switching network topology, which remains a challenging problem due to the complexity in the analysis of switching systems.

\subsubsection{Lyapunov Function Approach}

Although matrix theory is a relatively simple approach for stability analysis of the formation producing problem, it is not applicable in many formation producing scenarios, especially with nonlinear systems. It is then natural to consider the Lyapunov function approach.

By using the Lyapunov function approach, several typical formation producing scenarios have been studied, including the inverse agreement problem~\cite{DimarogonasKyriakopoulos09}, leaderless flocking and stabilization~\cite{Saber06,Do07,CuckerSmale07,TannerJadbabaiePappas07,BarnesFieldsValavanis09,KumarGargKumar10,ParkKimHa10,HaHaKim10,CuckerDong10,CuckerDong11,ZhangZhaiChen11}, and circular formation alike~\cite{SepulchrePaleyLeonard07,SepulchrePaleyLeonard08,ZhangLeonard07,Paley09,MellishNaporaPaley11,HernandezPaley10}. In the inverse agreement problem~\cite{DimarogonasKyriakopoulos09}, the objective is to force a team of agents to disperse in space. Roughly speaking, for the single-integrator kinematics~\eqref{eq:kinmatics}, the corresponding control input takes the form of
$
u_i(t)=\sum_{j=1}^n b_{ij}(\norm{x_i-x_j})[x_i(t)-x_j(t)],
$
where $b_{ij}(\cdot)$ is a nonnegative function. Assuming that each agent can communicate with all other agents within a radius $R$, the agents will disperse in space with the relative distance between any two agents being larger than $R$.

For the case of leaderless flocking, research has been conducted to stabilize a group of agents towards some desired geometric formation, where the inter-agent interaction is described directly or indirectly by some nonnegative potential function $V_{ij}(\norm{x_i-x_j})$ regardless of the final group velocity. Some notable properties of $V_{ij}(\norm{x_i-x_j})$ includes: (i) $V_{ij}(\norm{x_i-x_j})$ achieves its minimum when $\norm{x_i-x_j}$ is equal to the desired inter-agent distance between agents $i$ and $j$, (ii) $V_{ij}(\norm{x_i-x_j})$ increases as $\norm{x_i-x_j}$ decreases from the desired distance to zero and $V_{ij}(\norm{x_i-x_j})$ could approach infinity as $\norm{x_i-x_j}$ approaches zero, and (iii) $V_{ij}(\norm{x_i-x_j})$ increases as $\norm{x_i-x_j}$ increases from the desired distance to the maximum communication range $R$. The basic idea behind the potential function $V_{ij}(\norm{x_i-x_j})$ is to drive the inter-agent distance to the desired value while avoiding possible inter-agent collision. The corresponding control law for each agent is usually chosen to be the same as the corresponding consensus algorithm except that the $x_i-x_j$ term is replaced by $\nabla_{x_i} V(\norm{x_i-x_j})$ here. A fundamental limitation is that all agents will normally converge to some (constant) inter-agent configuration locally in the sense that some nonnegative potential function achieves its local minimum. Accordingly, the inter-agent distance might not converge to the desired value globally. It is an interesting future research topic to study how to ensure the desired inter-agent distance be achieved globally under a properly designed control algorithm. In addition, the network topology associated with a team of agents is usually assumed to be undirected, which is not applicable to many practical systems which are described by directed networks.

For the case of circular formation and the like, the main research in~\cite{SepulchrePaleyLeonard07,SepulchrePaleyLeonard08,PaleyLeonardSepulchre08,ZhangLeonard07,MellishNaporaPaley11} was devoted to the collective motion for nonhonolomic mobile robots with dynamics given in~\eqref{eq:nonholonomic}. Denote $r_i= x_i+\iota y_i$, where $\iota=\sqrt{-1}$. Then,~\eqref{eq:nonholonomic} becomes
$
\dot{r}_i = u_ie^{\iota \theta_i},~ \dot{\theta}_i=\omega_i,~ i=1,\cdots,N.
$
Due to the nature of the nonlinear dynamics, a consensus-like algorithm often renders a circular-like ultimate formation where the trajectories of all agents are circular and the relative phase difference (namely, $\theta_i-\theta_j$) is constant. The current work mainly focuses on the case when all agents share a common unit speed. Similar circular-like formation was analyzed in~\cite{Paley09,HernandezPaley10}, where the system dynamics are different from~\eqref{eq:nonholonomic} but share a similar nonlinearity. Due to the nonlinearity of the system dynamics, it is a challenging task to incorporate time delay, disturbances, quantization, etc, into the existing research.



\subsubsection{Graph Rigidity}

For a network with a given number of agents, the edges are closely related to the shape of formation. Roughly speaking, if enough information regarding edge distances for a team of agents is available, the geometric structure of all agents is determined. Then the graph for the agents is \emph{rigid}. According to~\cite{Laman70}, a graph of $n$ agents is rigid if at least $2n-3$ edge distances are available.


Motivated by the graph rigidity, research has been conducted in~\cite{SaberMurray02,HendrickxADB07,YuADF09,FidanHendrickxAnderson10,CaoYuAnderson11} to drive a group of agents to the desired inter-agent configuration by ensuring that a certain number of edge distances are identical to the desired ones. The graph rigidity recovery after loss of an agent has also been investigated. Compared with other formation producing algorithms which require edge vector information (\textit{i.e.,} $x_i-x_j$), less information is required in edge distance information (\textit{i.e.,} $\norm{x_i-x_j}$). As a tradeoff, some unstable equilibria, such as a collinear initial configuration (\textit{the initial states of all agents are linearly dependent}) and a common initial state (\textit{the initial states of all agents are identical}), rather than the desired inter-agent configuration, might exist. In practical applications, it is important to design proper control algorithms such that a team of agents can avoid converging to the unstable equilibria.


\subsubsection{Receding Horizon Approach}

Receding horizon control (RHC), known also as model predictive control, has been introduced into the formation stabilization problem. By nature, RHC is a finite-horizon optimization problem. The employment of RHC in the formation stabilization problem is motivated by the fact that RHC is more capable of dealing with constraints.

The main research on this topic~\cite{DunbarMurray06,Dunbar07,FrancoMPPR08} has been devoted to deriving proper distributed control algorithms for a team of agents such that they can reach some desired formation by optimizing some finite-horizon cost functions in scenarios with or without time delay. Because RHC is essentially an optimization-based control strategy, the distributed control algorithms are typically given by solving optimization problems. Therefore, more computational time is required by RHC than other control approaches. Therefore, the potential computation-induced time delay needs to be taken into consideration in practical applications.

\subsection{Formation Tracking}

Although formation control without a group reference is interesting in theory, it is more realistic to study formation control in the presence of a group reference because it may represent a control objective or a common interest of the whole group. This scenario is now reviewed in this subsection.

\subsubsection{Matrix Theory Approach}
Similarly to the case of formation producing, matrix theory is often used in the study of the formation tracking problem.

An interesting problem in formation tracking is to design a distributed control algorithm to drive a team of agents to track some desired state. For example, given the single-integrator kinematics, control algorithms were designed in~\cite{Ren07_SCL,CaoRenLi09}, where the algorithms are similar to those consensus algorithms except for that an extra term is introduced here due to the existence of the group reference. If properly designed, all agents can track the group reference accurately as reported in~\cite{Ren07_SCL}, with bounded tracking errors analyzed in~\cite{CaoRenLi09}, where a discretized version in~\cite{Ren07_SCL} was considered. It is worth mentioning that the group reference can be arbitrarily chosen as long as its derivative is bounded. In~\cite{XiangWeiLi09,WangHHJ10,PorfiriRobersonStilwell07}, the synchronization of a group of linear systems to the output of another linear exosystem was investigated with or without parameter uncertainties. In~\cite{Ren07_SCL,CaoRenLi09}, a general group reference was discussed while in~\cite{XiangWeiLi09,WangHHJ10,PorfiriRobersonStilwell07} a general system model was considered. How to solve formation tracking for general linear systems with a general group reference is still an open problem.

The formation tracking problem can be converted to a traditional stability problem by redefining the variables as the errors between each agent's state and the group reference. Then, the formation tracking problem is solved if the corresponding errors can be driven to zero. However, the formation producing problem, in general, cannot be solved in this way. Therefore, under a switching network topology, the formation tracking problem is generally easier to manager than the formation producing problem.


%


\subsubsection{Lyapunov Function Approach}

Due to the broad applications of the Lyapunov function approach in stability analysis, it has become an important tool in the study of the formation tracking problem as well.

Flocking with a dynamic group reference has been studied recently~\cite{SuWangLin09,ShiWangChu09,YuChenCao10_SCL,CaoRen11_TAC,WenDLC11}, where the objective is to design distributed control algorithms such that a team of agents move cohesively along the group reference. Compared with leaderless flocking where no specific final group velocity is required, the study of flocking with a dynamic group reference is much more challenging both theoretically and technically due to the existence of the dynamic group reference and the requirement on the cohesive movement of the agents along the dynamic group reference. In other words, the agents not only have to maintain some desired geometric formation but also need to follow the group reference as a whole. The combination of the two objectives makes the problem much more difficult than the leaderless flocking problem where only the first objective is involved. If enough information of the group reference is known, such as the acceleration and/or velocity information of the group reference, flocking with a dynamic group reference can be solved by employing a gradient-based control law~\cite{SuWangLin09,ShiWangChu09,YuChenCao10_SCL}. Another approach was proposed in~\cite{CaoRen11_TAC}, where a variable structure-based control law was used to solve the problem with less information required. Similarly to the study of the leaderless flocking problem, the existing approaches on flocking with a dynamic group reference can only reach a local minimization of certain potential functions because the potential function is generally unspecified but satisfies the conditions stated in Subsection~\ref{sec:formatioProducing}. Accordingly, the inter-agent distance is not identical to the desired one. However, the potential-based control can be easily designed to guarantee collision avoidance and maintain the initial inter-agent communication patterns. Nevertheless, it is still an open problem to accomplish the task with global inter-agent distance stabilization, collision avoidance, and initial communication pattern maintenance.

Formation control with a group reference was studied in both linear systems~\cite{NiCheng10,Ren10_TCST1,HuFeng10} and nonlinear systems~\cite{WangSlotine06,SongCaoYu10,Do08,GustaviHu08,DongFarrell08,DongFarrell09,Dong11,SedaTEMSS10,Ren10_TCST2} when the potential function $V(\norm{x_i-x_j})$ is replaced by some known functions, generally in the form of $\norm{x_i-x_j-d_{ij}}^2$, where $d_{ij}$ denotes the desired distance between agents $i$ and $j$. Briefly, the nonlinear systems studied in this case include nonholonomic mobile robots (see~\eqref{eq:nonholonomic})~\cite{Do08,DongFarrell08,GustaviHu08,DongFarrell09,Dong11}, rigid bodies (see~\eqref{eq:rigid-body})~\cite{SedaTEMSS10,Ren10_TCST2}, and linear systems with other nonlinear terms~\cite{WangSlotine06,SongCaoYu10}. Compared with the flocking problem, the problem studied here is relatively easier due to the known $V(\norm{x_i-x_j})$. In general, the inter-agent distance can be driven to the desired one. As a tradeoff, the collision avoidance and initial communication pattern maintenance need to be considered separately.

%
%


\subsection{Connectivity Maintenance for Consensus and Formation Control}

In both consensus and formation control problems, it is often assumed that the network topology satisfies certain fundamental
conditions, for example, is connected or has a directed spanning tree. However, a practical communication model is typically distance-based, i.e., two agents can communicate with each other if and only if their distance is smaller than a certain threshold, called \textit{communication range}. This is particularly true for sensor networks. In order to guarantee consensus or formation control be achieved asymptotically, a connectivity maintenance mechanism is essential, which has been studied recently.

The main approach to maintaining the connectivity of a team of agents is to define some artificial potentials (between any pair of agents) in a proper way such that if two agents are neighbors initially then they will always communicate with each other thereafter~\cite{SuWangChen10,ZavlanosTJP09,AjorlouMomeniAghdam10,GustaviDEH10,SuWangLin09_automatica,SuWangChen09,JiEgerstedt07, ZavlanosPappas08,ShuckerMurpheyBennett08,ZavlanosPappas07,CaoRen11_TAC}. In general, the artificial potential between a pair of agents grows to be sufficiently large (could be unbounded) when the distance between them increases to be equal to the communication range. For properly designed control algorithms, which are usually composed of the gradients of the artificial potentials, the total artificial potential is nonincreasing. This then indicates that the initial communication patterns can be preserved because otherwise the total potential will become larger than the initial total artificial potential, as soon as some communication pattern is broken. Although this approach provides a systematic way to guarantee the connectivity, the corresponding control algorithms might require infinite large control inputs, which is not practical. Meanwhile, it is not even necessary to always maintain the initial communication patterns in order to guarantee the connectivity. Therefore, how to find a more effective way to guarantee connectivity deserves further investigation. In contrast to the studies in~\cite{SuWangChen10,ZavlanosTJP09,AjorlouMomeniAghdam10,GustaviDEH10,SuWangLin09_automatica,SuWangChen09,JiEgerstedt07, ZavlanosPappas08,ShuckerMurpheyBennett08,ZavlanosPappas07,CaoRen11_TAC}, the authors in~\cite{SunHuang09} investigated an interesting problem where the number of initially existing communication patterns plays a role in the connectivity maintenance for the consensus problem with single-integrator kinematics~\eqref{eq:kinmatics} and control input~\eqref{eq:control}. Roughly speaking, if the initial graph is ``sufficiently'' connected in the sense that each agent has at least a certain number of neighbors, consensus can be guaranteed to be achieved. Note that the result can only be applied to systems with single-integrator kinematics therefore further investigation is expected for systems with high-order linear dynamics or nonlinear dynamics.

\subsection{Remarks}
Current research on formation control mainly focuses on a fixed formation where the inter-agent distance is fixed. Considering practical applications, however, it might require the formation be adaptive with respect to the events performed by the team of agents. In addition, it is important to consider constraints, such as input saturation, quantization, and power limitation, in the formation control problem. Meanwhile, the robustness is another important factor that deserves considerable attention in real applications where noise and disturbances exist.

In terms of connectivity maintenance for consensus and formation control, research has been devoted to continuous-time systems but practical systems are more suitable to be modeled in a discrete-time setting, which makes the study of connectivity maintenance for discrete-time systems more meaningful. In general, the connectivity maintenance for discrete-time systems is more challenging due to the fundamental limitation of the corresponding control input, which is usually piecewise constant rather than continuous.

\section{Optimization} \label{sec:optimization}

Optimization is an important subject in the studies of control systems. The main objective of optimization is to find an optimal strategy subject to some given constraints such as cost functions.
Recently optimization in distributed multi-agent coordination has been studied concerning convergence speed and some specific cost functions, which are respectively reviewed below.

\subsection{Convergence Speed}

As discussed above, one important problem in consensus is the convergence speed, which characterizes how fast consensus can be achieved therefore is desirable to optimize. In this regard, the convergence speed is the cost function to be maximized.

Consider a group of $n$ agents with dynamics described by the single-integrator kinematics~\eqref{eq:kinmatics}. Equipped with~\eqref{eq:control}, the dynamical system~\eqref{eq:kinmatics} can be written in a matrix form, as
\begin{equation}\label{eq:clp}
  \dot{X}(t)=-\Lcal X(t),
\end{equation}
where $X(t)=[x_1(t),\cdots,x_n(t)]^T$ and $\Lcal$ is the Laplacian matrix. For a network with a fixed topology, $\Lcal$ is a constant matrix.


Motivated by the observation that the smallest nonzero eigenvalue of the Laplacian matrix, $\lambda_2(\Lcal)$, determines the worst-case convergence speed~\cite{SaberMurray04}, research has been conducted to maximize the convergence speed~\cite{KimMesbahi06,Kim10} by choosing optimal weights associated with edges. In contrast to~\cite{KimMesbahi06,Kim10}, where the systems are assumed to have single-integrator kinematics, optimization of the convergence speed for double-integrator dynamics was considered in~\cite{CarliCSZ11}, where the convergence speed is defined in a similar way to the $\lambda_2(\Lcal)$ for the case with single-integrator kinematics. It is worth mentioning that optimal convergence for general linear and nonlinear systems is still an open problem.

Other than $\lambda_2(\Lcal)$, another commonly used measure for the convergence speed is given by~\eqref{eq:convspeed}. The corresponding optimization problem is
\begin{align}\label{eq:opti}
  \max_{u_i\in\Ucal}~\rho,
\end{align}
where $\rho$ is defined in~\eqref{eq:convspeed} and $\Ucal$ is a set of admissible controllers. Existing research in~\cite{XiaoBoyd04,KimGuPostlethwaite09} focuses on the case when all agents converge to the average of the initial states, \textit{i.e.,} $X^\star = [\frac{1}{n}\sum_{i=1}^n x_i(0)]\onebf$. For an arbitrary fixed or switching network topology, the optimization problem~\eqref{eq:opti} is challenging if $X^\star$ is unknown. But if $X^\star$ is chosen as $ [\frac{1}{n}\sum_{i=1}^n x_i(t)]\onebf$, then the problem becomes much easier.


\subsection{Specific Cost Functions}

In addition to the fastest convergence speed requirement, various cost functions are also subject to minimization.

One interesting problem studied in this setting is distributed multi-agent optimization, which is motivated by solving one challenge in wireless sensor networks, namely, to estimate the environment parameters and/or some functions of interest, such as temperature and source location~\cite{RabbatNowak04}. As a simple strategy, each sensor node can send its data to an existing central location which can then process the data if it is sufficiently powered. However, due to the limited power resources and communication capabilities, this strategy is often not applicable. An alternative approach to achieving a similar objective is to estimate the environment parameters and/or some functions of interest locally, which requires much less communication bandwidth and power. In wireless sensor networks~\cite{RabbatNowak04}, the estimation problem is usually modeled as a distributed multi-agent optimization problem. Roughly speaking, the objective of distributed multi-agent optimization is to cooperatively minimize the cost function
$\sum_{i=1}^{n}f_{i}(x),$
where the function $f_{i}:\re^{n}\mapsto\re$ represents the cost
of agent $i$, known by this agent only, and $x\in\re^{m}$ is a decision vector. In~\cite{RabbatNowak04}, an incremental subgradient approach was used to solve the optimization problem for a ring type of network. It should be noted that~\cite{RabbatNowak04} does not provide much discussion on the optimization problem under other types of network topologies.

Ref.~\cite{NedicOzdaglar09} was probably the first paper studying the distributed multi-agent optimization problem under a consensus-based framework. The problem considered therein is formulated so as to
$$
\begin{array}{lll}
\mathrm{minimize}    &  & \sum_{i=1}^{n} f_{i}(x),                                                 \\
\mathrm{s.t.}&  & x \in \re^{n},
\end{array}
\label{mao}
$$
where each $f_{i}:\re^{n}\rightarrow\re$ is assumed to be a convex function. Inspired by the average consensus algorithm and the standard subgradient method, a consensus-like algorithm was proposed as
\begin{equation}\label{eq:opti-sub}
x_{i}(k+1) = \sum_{j=1}^{n} a_{ij}(k) x_{j}(k) - \alpha g_{i}(x_{i}(k)),
\end{equation}
where $\alpha$ is the step size and $g_{i}(x_{i}(k))$ is the subgradient of $f_i(x)$ at $x=x_i(k)$. In~\cite{RabbatNowak04}, $\sum_{j=1}^{n} a_{ij}(k) x_{j}(k)$ in~\eqref{eq:opti-sub} was replaced by $x_{i-1}(k)$ with $x_0(k)=x_n(k-1)$. Note that the algorithm~\eqref{eq:opti-sub} can only find sub-optimal solutions, determined by the constant
step size $\alpha$. Further results on this topic can be found in~\cite{NedicOzdaglarParrilo10,LobelOzdaglar11,Nedic11}, where a similar distributed multi-agent optimization problem was studied within various scenarios, under constraints~\cite{NedicOzdaglarParrilo10}, over random networks~\cite{LobelOzdaglar11}, and with broadcast-based communications in an asynchronous setting~\cite{Nedic11}. In the existing literature, time delay and disturbances have not been taken into consideration. Therefore, it is important to consider time delay and disturbances in the distributed multi-agent optimization problem due to their wide existence in wireless sensor networks.

In addition to the distributed multi-agent optimization problem, where the cost function is a sum of a series of convex functions, distributed optimization has also been considered for both infinite-horizon cost functions~\cite{CaoRen09_SMCB,BorrelliKeviczky08,KazerooniKhorasani10,WangXin10,Dong10} given by
$
J_i = \int_{0}^\infty [X^T(t)QX(t)+U^T(t)RU(t)]\text{d}t
$
and finite-horizon cost functions~\cite{HuPrandiniTomlin06,FrancoParisiniPolycarpou07,TrecateGMS09,StewartVRWP10} given by
$
J_f = \int_{0}^{t_f} [X^T(t)QX(t)+U^T(t)RU(t)]\text{d}t,
$
where $X\in\re^n$ is the state, $U\in\re^n$ is the control input, and $t_f$ is a positive constant. It is worth mentioning that the RHC approach discussed in Section~\ref{sec:formatioProducing} typically has finite-horizon cost functions. Differing from the research reported in~\cite{NedicOzdaglar09,NedicOzdaglarParrilo10,LobelOzdaglar11,Nedic11}, which is to find the optimal estimated state, the objective here is to find the optimal control laws subject to the minimization of certain cost functions. Due to requirements of optimizing the cost functions when designing the control laws, the computational complexity becomes an important problem to study. Meanwhile, the network topology plays a significant role in the optimization problem with certain cost functions, therefore it is also important to optimize the network topology subject to certain cost functions.



\section{Distributed Task Assignment} \label{sec:task}

Distributed task assignment refers to the study of task assignment of a group of dynamical agents in a distributed manner, which can be roughly categorized into coverage control, scheduling, and surveillance. Compared with the previous studies discussed in Sections~\ref{sec:consensus},~\ref{sec:formation}, and~\ref{sec:optimization}, distributed task assignment focuses on the three task-oriented research problems, where each topic has its unique features.

\subsection{Coverage Control}

Coverage control is an active research direction in mobile sensor networks, where the objective is to properly assign the mobile sensors' motion in order to maximize the detection probability.

Let $Q$ be a convex space with $\phi$ represent the distribution density function which indicates the probability that some event takes place over $Q$~\cite{CortesMKB04_TRA}. Consider a group of $n$ mobile sensors whose locations are specified by $P=[p_1,\cdots,p_n]$. The sensor performance at a point $q$ degrades with respect to the distances $\norm{q-p_i}$, which are all described by a nondecreasing differentiable function, $f$. The coverage control problem is to find a local controller for each mobile sensor such that the following cost function is minimized:
$
  J = \sum_{i=1}^n\int f(\norm{q-p_i})\phi(q)dq\,.
$
This coverage control by nature is an optimization problem. Main research on this topic was reported in, e.g.,~\cite{CortesMartinezBullo05,GaoCortesBullo08,LaventallaCortes09,KwokMartinez11,WangHussein10,ChoiHorowitz10,Cortes10,KwokMartinez10}, where the coverage control problem was considered in two directions, namely, analysis of coverage control under various practical constraints, such as limited sensing/communication capacities~\cite{CortesMartinezBullo05}, load balancing~\cite{Cortes10}, and nonholonomic mobile robots~\cite{KwokMartinez10}, and algorithms for coverage control~\cite{KwokMartinez11}. Noting that time delay and uncertainties have not been considered in the coverage control problem, it is interesting to consider the effect of time delay and uncertainties in the coverage control problem. Moreover, the density function $\phi$ might be time-dependent in real systems, which is another interesting research topic for further study.



\subsection{Scheduling}

Another interesting topic in distributed task assignment is distributed scheduling, which refers to the scheduling of a group of dynamical agents in a distributed manner. Distributed scheduling has many potential applications in military and civilian sectors, and can be roughly categorized into two typical problems, namely sequence optimization~\cite{JinShimaSchumacher06} and task
allocation~\cite{FinkePassinoSparks06,MooreLucarelli07,ChoiBrunetHow09,ReveliotisRoszkowska11,SayyaadiMoarref11}. The objective of sequence optimization is to schedule a team of agents such that some metrics can be optimized. For instance, in~\cite{JinShimaSchumacher06}, an optimal scheduling sequence was designed to fuel a group of UAVs via dynamic programming, where the metric is the total spending time. The objective of task allocation is to distribute certain number of tasks to a team of agents such that they can balance the total tasks. If the number of tasks for each agent is considered a variable in the consensus problem, the task allocation might be viewed as a consensus problem except that a limitation on the total number of tasks for all agents exists. A notable feature of the distributed task assignment problem is that various constraints may exist due to the physical properties associated with the agents. In view of the difference objectives for sequence optimization and task allocation, it is an interesting topic to consider combining both objectives of the two problems simultaneously.


\subsection{Surveillance}

Distributed surveillance means to monitor a certain area by using a group of mobile agents in distributed coordination. Distributed surveillance has a number of potential applications, such as board security guarding, forest fire monitoring, and oil spill patrolling.

The main motivation of distributed surveillance is that a team of agents can monitor a given (large) area more effectively than a single agent when the team of agents works in a cooperative fashion. Accordingly, an important research problem in distributed surveillance is to design environment-based cooperative control laws for all coordinated agents such that the given area can be monitored efficiently. Recent research on distributed surveillance has been reported in~\cite{AnisiOgrenHu10,KingstonBeardHolt08,MoorePassino08,GilPassinoCruz08,PeasgoodClarkPhee08,YangPolycarpouMinai07,FinkePassino07,TangOzguner08}, where a number of physical limitations were identified and considered such as time delay and uncertainties~\cite{GilPassinoCruz08,FinkePassino07}, collision avoidance between agents~\cite{PeasgoodClarkPhee08}, and heterogeneously distributed agents~\cite{FinkePassino07}. The current research is conducted under the assumption that each agent has enough power such that any designed control law can be applied. However, due to the power constraints, each agent might be subjected to constraints such as bounded control input, limited distance to travel, and finite accuracy level, etc., therefore it is interesting and important to consider these limitations in the distributed surveillance problem in the future.

\section{Estimation} \label{sec:est}

Due to the absence of global information, used for achieving group coordination in many cases, a distributed estimation scheme is often needed for estimation.

The first problem is to design local distributed estimators such that some global information can be estimated asymptotically or in finite time. The second problem is to design local controllers based on the local estimators such that the closed-loop system is stable. The estimation-based distributed control is essentially a combination of both centralized control and distributed control in such a way that distributed control is used in the estimation of some global information and the centralized control idea is used for local controllers design. The estimation-based distributed control strategy often inherits the merits of both centralized control and distributed control. However, it is worth emphasizing that a closed-loop system with distributed estimators is much more complicated to design than one without distributed estimators.

Main research on this topic has been reported in~\cite{CaoRenMeng10,ZhangLeonard10,YangFGLSS10,ChoiOhHorowitz09,LynchSYF08,YangFreemanLynch08,SmithHadaegh07}, where the joint estimation and control problem was considered subject to disturbances~\cite{ZhangLeonard10,ChoiOhHorowitz09,LynchSYF08} or without disturbances~\cite{CaoRenMeng10,YangFGLSS10,YangFreemanLynch08,SmithHadaegh07}. In~\cite{CaoRenMeng10,ZhangLeonard10,YangFGLSS10,ChoiOhHorowitz09,LynchSYF08,YangFreemanLynch08,SmithHadaegh07}, a joint estimation and control problem is solved in the sense that the distributed estimator is used in the design of proper control algorithms such that certain global objective can be achieved. Without the aid of distributed estimators, the control design is very hard and even impossible. As can be noticed from~\cite{StankovicStankovicStipanovic09,CalafioreAbrate09,DemetriouHussein09,MarechalGorcePierrot10,StankovicStankovicStipanovic11,Yu10,SubbotinSmith09}, the distributed estimation problem has been considered without much discussion on specific control problems. In general, the joint estimation and control idea has provided an important approach in the study of distributed multi-agent coordination where only neighbor-based information is not sufficient for the controllers design. On the other hand, in real applications properly designed distributed estimators might be used to replace some expensive sensors.

In general, it remains a challenging problem to study task-oriented coordination control systems where the use of distributed estimation is either necessary or an appropriate replacement of certain expensive measurement devices, at the costs of difficult control system design and complex system stability analysis. Moreover, physical limitations such as bounded control input, asynchronous communication, and information quantization, could potentially reduce the applicability of the joint estimation and control scheme in various distributed multi-agent coordination systems.

\section{Discussions} \label{sec:conclusion}

This article has reviewed some recent research and development in distributed multi-agent coordination. In addition to the theoretical results reviewed above, many experiments have also been conducted to validate the theoretical designs and analysis, as can be found in~\cite{HiraiMiyazaki06,KunwarBenhabib06, AntonelliChiaverini06,BaldassarreTBMDN07,RenCBSC08,GradyChristensenDorigo09},
to name just a couple of representative reports. Although the existing theoretical research and experiments have solved a number of technical problems in distributed multi-agent coordination, there are still many interesting, important and yet challenging research problems deserving further investigation. Some of them are briefly summarized as follows:
\begin{itemize}
\item {\it Quantization effects in distributed coordination algorithms.} The current research efforts focus on studying distributed coordination problems with control inputs and measurements being analog signals with continuous values. However, digital signal processing techniques require digital inputs and sampled-data measurements. Although quantization effects have been studied in several coordination problems, the quantization effect in some other distributed coordination problems remain unsolved and even untouched.
\item {\it Optimization with integrated individual and global cost functions.} Optimization problem in distributed coordination has been studied with various cost functions. In real systems, each individual agent has both local and global objectives, contributing to an integration of both individual and global cost functions. Therefore, optimizing a combined objective is more realistic but also more challenging. Another interesting problem is to investigate how to balance the individual cost functions and the global cost function toward a common objective.
\item {\it Intelligent coordination.} Intelligent coordination refers to the distributed coordination of a team of agents in the presence of artificial intelligence, namely, each agent is intelligent in the sense that they can choose the best
    possible responses based on its own objective. Intelligent coordination has potential applications not only in engineering and technology but also in economics and social studies. Although research problems, such as pursuer-invader problem~\cite{CaoTLGW06,BopardikarBulloHespanha09,BopardikarBulloHespanha08,KollingCarpin10} and the game theory in distributed coordination~\cite{LiangXiao10,Gu08,BausoGiarrePesenti09,BausoGiarrePesenti08,VrancxVerbeeckNowe08,HuangCainesMalhame07}, have been studied recently, there are still many open questions especially the understanding of group behaviors in the presence of agent intelligence. One interesting problem is how to interpret the underlying complex networks as well as to stabilize and optimize the network in the presence of agent intelligence.
\item {\it Competition and cooperation.} Today, most research is conducted based on local cooperation but not competition. This posts an obvious limitation because competition not only exists but also plays an positive role in group coordination. For example, due to the lack of competition, the final states of the traditional consensus algorithms are always limited to be within some region in the state space determined by the initial agent states. One interesting question is how to introduce competition into distributed coordination so as to arrive at different desired regions and to improve the system performance by rewarding different agents with different benefits.
\item {\it Centralization and decentralization.} Although decentralization shows obvious advantages over centralization, such as scalability and robustness, decentralization also has its own drawbacks. One shortcoming is that, under decentralized protocols, some agents cannot predict the group behavior based only on the available local information. Consequently, some group behavior cannot be controlled. As a sensible example, current economic crisis actually illustrates some disadvantages of behavioral decentralization. One interesting question, therefore, is how to balance decentralization and centralization so as to further improve the overall systems performance.
\end{itemize}

\bibliographystyle{IEEEtran}
\bibliography{survey-refs-arxiv}

\end{document}